\documentclass[twoside,reqno,11pt]{fcaa-var} %

\usepackage{graphicx}
\usepackage{epsfig}
\usepackage{amsthm}
\usepackage{amsmath}
\usepackage{latexsym}
\usepackage{amsfonts}
\usepackage{amssymb}

\newtheoremstyle{theorem}
  {15pt}          
  {15pt}  
  {\sl}  
  {\parindent}
  {\sc}  
  {. }   
  { }    
  {}     
\theoremstyle{theorem}
\newtheorem{lemma}{Lemma}[section]
\newtheorem{theorem}{Theorem}[section]

\newtheoremstyle{defi}
  {15pt}          
  {15pt}  
  {\rm}  
  {\parindent}     
  {\sc}  
  {. }    
  { }    
  {}     
\theoremstyle{defi}
\newtheorem{definition}{Definition}[section]

 \def\proofname{\indent\rm P\ r\ o\ o\ f.}
 \def\proofend{\hfill$\Box$}
 \def\qed{\hfill$\Box$} 

 \usepackage{hyperref} 
\usepackage{color}
 \def\theequation{\arabic{section}.\arabic{equation}}

 \def\themycopyrightfootnote{\vspace*{3pt}
 \ \, \,  
   \par  \noindent 
   \hfill  \vspace*{-36pt}
  }

 \title[FDE with dist.-order material derivative \dots]{Fractional diffusion equation with distributed-order material derivative. Stochastic foundations. \\ [3pt] }
 \author[\normalsize M. Magdziarz, M. Teuerle]{\normalsize Marcin Magdziarz $^1$, Marek Teuerle $^2$}

\newcommand{\E}{\mathrm{{E}}}

\newcommand{\Prob}{\mathrm{{P}}}

\newcommand{\dd}{\mathrm{d}}
\newcommand*{\eqb}{\begin{eqnarray}}
\newcommand*{\eqe}{\end{eqnarray}}
\newcommand*{\al}{\alpha}

\newcommand{\fconv}{\Rightarrow}

\newcommand{\indyk}{{\mathbf{1}}}

\begin{document}

 \vbox to 2.5cm { \vfill }


 \bigskip \medskip

 \begin{abstract}

In this paper we present stochastic foundations of fractional dynamics driven by fractional material derivative of distributed order-type. Before stating our main result we present the stochastic scenario which underlies the dynamics given by fractional material derivative. Then we introduce a Levy walk process of distributed-order type to establish our main result, which is the scaling limit of the considered process. It appears that the probability density function of the scaling limit process fulfills, in a weak sense, the fractional diffusion equation with material derivative of distributed-order type.
 \medskip

 \smallskip


 \end{abstract}

 \maketitle

 \vspace*{-16pt}



\section{Introduction}\label{sec:1}
\setcounter{section}{1}
\setcounter{equation}{0}\setcounter{theorem}{0}

The scaling limit of the standard random walk with jumps having finite variance
leads to the classical diffusion equation
\[
\frac{\partial p}{\partial t}=\frac{1}{2}\frac{\partial^2 p}{\partial x^2} .
\]
By adding the heavy-tailed waiting times which are independent of the jumps, we obtain the continuous-time random walk (CTRW), whose long-time limit is governed by the fractional diffusion equation \cite{Metzler_Klafter}
\[
\frac{\partial p}{\partial t}=
\frac{1}{2} {_0D_t^{1-\al}} \frac{\partial^2 p}{\partial x^2}  .
\]
The operator
$_0D_t^{1-\al}$, $0<\al<1, \;f\in C^1([0,\infty)),$ is the fractional derivative of the Riemann-Liouville type \cite{samko}.

On the other hand, if both jumps and waiting times are independent and chosen from the power-law distribution,
the density of the corresponding scaling limit solves the space-time fractional diffusion equation \cite{Metzler_Klafter}.
The main difficulty underlying such models with heavy-tailed jumps is the infinite second moment.
As a consequence, the mean square displacement is not well defined.
To overcome this difficulty the so-called Le\'vy walks were introduced \cite{Klafter1,Klafter2} (see also \cite{LW_review}
for a recent review).
L\'evy walk is a special type of CTRW, for which the length of each jump is equal to
the length of the preceding waiting time. This
strong spatiotemporal coupling assures that the L\'evy walk has
finite moments of all orders (long jumps are
penalized by long waiting times preceding the jump).
Thus, even if the jumps are heavy-tailed,
the mean square displacement is well-defined.
Moreover, the dynamics of the L\'evy walk densities is
governed by the so-called fractional material derivatives $\left( \frac{\partial}{\partial t} \mp
\frac{\partial}{\partial y} \right)^\al$, see \cite{Sokolov_Metzler}. In the
Fourier-Laplace space they are given by
\[
\mathcal{F}_y\mathcal{L}_t\left\{\left( \frac{\partial}{\partial t} \mp \frac{\partial}{\partial y} \right)^\al f(y,t)\right\}=(s\mp ik)^\al f(k,s).
\]

In this paper we introduce the L\'evy walk model with waiting times
displaying very slow logarithmic decay of the tails. This type of
distributions leads to the so-called ultraslow diffusion, which was analyzed in \cite{Chechkin,Meerschaert_ultraslow}.
Because of the strong coupling, in spite of the fact that the jumps can be extremely long, the introduced model
has finite moments of all orders. Moreover, it explains the stochastic origins of fractional dynamics driven by fractional material derivative of distributed order-type.

The paper is organized as follows: in the next section we recall the notion of fractional material
derivative and its relation with the classical L\'evy walks.
Next we define he L\'evy walk model with waiting times
having logarithmic decay of the tails. We derive the corresponding scaling limit
in the cases of first-wait and first-jump scenarios.
Finally, we introduce a pseudo-differential operator called \emph{fractional material derivative of distributed-order type}
and show that the density of the obtained diffusion limit solves the corresponding
fractional diffusion equation.

\section{Stochastic foundations of fractional dynamics with fractional material derivative}\label{sec:2}
\setcounter{section}{2}
\setcounter{equation}{0}\setcounter{theorem}{0}
In the following section we recall the formal definition of multidimensional L\'evy walk (LW)  and the corresponding overshooting L\'evy walk (OLW). The scaling limits of these processes lead to the fractional dynamics which involves fractional material derivatives.

\subsection{Definitions of L\'evy walk and overshooting L\'evy walk}

LW process arises as a special case of CTRW. CTRW is a stochastic process, which generalizes the classical random walk \cite{Pearson, Feller} by assuming that the consecutive random jumps are separated by random waiting times \cite{Feller,Weiss,Montroll,Sokolov_book,Metzler_Klafter}.
\begin{definition}
\label{def:lwolw}
Let $\{(T_i,\mathbf{J}_i)\}_{i\geq 1}$ be a sequence of independent and identically distributed (IID) random vectors such that the following conditions hold:
\begin{itemize}
\item[a)] $\Prob(T_i>0)=1$ and $T_i$ is heavy-tailed with index $\alpha\in(0;1)$, i.e.
\begin{equation}
\Prob(T_i>t)\approx t^{-\alpha} \qquad\mathrm{as}\quad t\to\infty,
\label{heavytailed}
\end{equation}
\item[b)] $\Prob(\mathbf{J}_i\in\mathbb{R}^d)=1$ and for each $\mathbf{J}_i$ we have
\begin{equation}
\mathbf{J}_i=\mathbf{V}_iT_i,
\label{jumps}
\end{equation}
where $\{\mathbf{V}_i\}_{i\geq 1}$ is an IID sequence of non-degenerate unit vectors from $\mathbb{R}^d$, governing the direction of jumps. Sequences $\{\mathbf{V}_i\}_{i\geq 1}$ and $\{T_i\}_{i\geq 1}$ are assumed mutually independent.
\end{itemize}
Next, let $N_t=\max\{n:T_1+T_2+\ldots+T_n\leq t\}$ be the counting process corresponding to $\{T_i\}_{i\geq 1}$. Then, the processes
\begin{equation}
\mathbf{R}(t)=\sum_{i=0}^{N_t}\mathbf{J}_i,
\label{lw}
\end{equation}
\begin{equation}
\mathbf{\tilde{R}}(t)=\sum_{i=0}^{N_t+1}\mathbf{J}_i,
\label{olw}
\end{equation}
are called LW and OLW, respectively.
\end{definition}

Physically, the process defined in (\ref{lw}) has a very intuitive behaviour and can be used to model random phenomena at the microscopic level. Namely, a particle whose position is described by LW, starts its motion at the origin and stays there for the random time $T_1$, then it makes the first jump $\mathbf{J}_1$. Next, it stays at the new location for the random time $T_2$ before making next jump $\mathbf{J}_2$, and then the scheme repeats. Due to relation (\ref{jumps}), the length of each jump $\mathbf{J}_i$ is equal to the length of the preceding waiting time $\mathbf{T}_i$, while the direction of each jump $\mathbf{J}_i$ is governed by the random vector $\mathbf{V}_i$. Therefore the position of the particle $R_t$ at any moment $t$ satisfies $\|R_t\|\leq t $, where $\|\cdot\|$ is $d$-dimensional Euclidean norm.

As the reader can notice in the definition of OLW the counting process $N_t$ is replaced by $N_t+1$. This change has a important influence on the trajectories and properties of the underlying OLW process.
In contrary to LW process, the particle running OLW performs the first jump $\mathbf{J}_1$ immediately after it starts the motion, and waits at the new location for the time $T_1$, then it performs the second jump $\mathbf{J}_2$ to relocate to the next position and waits there for the time $T_2$, and then that scheme repeats. In fact the OLW process consists of 'jump-wait' events (also called 'first jump' events), while LW process consists of 'wait-jump' events ('first wait' events). One of the main differences between LW and OLW processes is that the second moment of the LW process is finite, while it is infinite for OLW process \cite{Magdziarz2015,Magdziarz2012,Teuerle2012,TeuerleOver2012}

\subsection{Asymptotic properties of L\'evy walks and overshooting L\'evy walks}

In this section we present the asymptotic behaviour of LW and OLW processes. We give the detailed description of the corresponding scaling limit processes.

Recall that $\{\mathbf{X}(t)\}_{t\geq 0}$ is a $d$-dimensional L\'evy process if its Fourier transform is given by (L\'evy-Khinchin representation, see \cite{sato}):
\begin{equation}
\label{levy-khinchin}
\E \exp\{i\langle\mathbf{k},\mathbf{X}(t)\rangle\}=\exp\{t\psi(\mathbf{k})\},\qquad \mathbf{k}\in\mathbb{R}^d
\end{equation}
with the characteristic exponent $\psi(\mathbf{k})$ given by:
\begin{equation}
\label{levy-khinchin_exponent}
\psi(\mathbf{k})=i\langle\mathbf{k},\mathbf{a}\rangle-\frac{1}{2}\langle\mathbf{k},\mathbf{Qk}\rangle+\int_{\mathbf{x}\neq\mathbf{0}}{\left(\exp{i\langle\mathbf{k},\mathbf{x}\rangle}-1-\frac{i\langle\mathbf{k},\mathbf{x}\rangle}{1+\|\mathbf{x}\|^2}\right)\nu(\dd x)}.
\end{equation}
Here, $\langle\cdot,\cdot\rangle$ is the scalar product in $\mathbb{R}^d$,  $\mathbf{a}\in\mathbb{R}^d$ is the drift parameter, $\mathbf{Q}$ is the symmetric and positive defined $d$-dimensional square matrix determining the Gaussian part of $X(t)$, and $\nu$ is the so-called L\'evy measure. It is the $\sigma$-finite Borel measure on $\mathbb{R}^d\setminus\{\mathbf{0}\}$ such that $\int_{\mathbf{x}\neq\mathbf{0}}\min(\|\mathbf{x}\|^2,1)\nu(\dd x)<\infty$.
The triplet $[\mathbf{a},\mathbf{Q},\mathbf{\nu}]$ , which is called the L\'evy triplet, uniquely determines the L\'evy process $X(t)$.

Let us introduce the following notation: denote by $\mathbf{X}^-(t) = \lim_{s\nearrow t}\mathbf{X}(s)$ the left-continuous version of a right-continuous processes $\mathbf{X}(t)$, and by
$\mathbf{Y}^+(t) = \lim_{s\searrow t}\mathbf{Y}(s)$ the right-continuous version of any left-continuous process $\mathbf{Y}(t)$.
\begin{theorem}
\label{th:fconv}
Let $\mathbf{R}(t)$ and $\mathbf{\tilde{R}}(t)$ be the LW and OLW processes from Def. \ref{def:lwolw}. The following convergences hold in $\mathbb{J}_1$ topology as $n\to\infty$:
\begin{eqnarray}
\label{th:fconv1}
n^{-1/\alpha}{\mathbf{R}\left(n^{1/\alpha} t\right)}\Rightarrow \left( \mathbf{L}^{-}_{\alpha}\left( S^{-1}_{\alpha}(t)\right)\right)^+,\\\nonumber n^{-1/\alpha}{\mathbf{\tilde{R}}\left(n^{1/\alpha} t\right)}\Rightarrow  \mathbf{L}^{}_{\alpha}\left( S^{-1}_{\alpha}(t)\right).
\end{eqnarray}
Here, the process $S^{-1}_{\alpha}(t)$ is defined as follows:
\begin{equation}
S^{-1}_{\alpha}(t)=\inf\{\tau \geq 0: S_{\alpha}(\tau)>t\}.
\label{inversestablesubordinator}
\end{equation}
The dependence structure between the L\'evy processes $\mathbf{L^{}}_{\alpha}(t)$ and $S^{}_{\alpha}(t)$ is given by their joint L\'evy triplet,
\begin{equation}
\label{jointlevytriplet1}
\left[\ \int\limits_{\mathbf{x}\neq \mathbf{0}}\frac{\mathbf{x}}{1+\|\mathbf{x}\|^2}\mathbf{\nu}_{(\mathbf{L}^{}_{\alpha},S^{}_\alpha)}(\dd \mathbf{x}),\mathbf{0},\mathbf{\nu}_{(\mathbf{L}^{}_{\alpha},S^{}_\alpha)}\right],
\end{equation}
with the following L\'evy measure
\begin{equation}
\label{jointlevymeasure1}
\mathbf{\nu}_{(\mathbf{L}^{}_{\alpha},S_\alpha)}(\dd \mathbf{x})=\int\limits_{\mathbb{S}^{d-1}}\delta_{\mathbf{u} t}(\dd \mathbf{x}_1)  \nu(\dd t)\mathbf{\Lambda}(\dd \mathbf{u}).
\end{equation}
Here $\mathbf{x}=(\mathbf{x}_1,t)\in \mathbb{R}^d\setminus\{\mathbf{0}\}\times(0,\infty)$, $\mathbf{u}=\mathbf{x}_1/\|\mathbf{x}_1\|\in \mathbb{S}^{d-1}$, $\mathbb{S}^{d-1}$ is the (d-1)-dimensional sphere, $\delta_{s}(\cdot)$ is a Dirac delta function at point s, $\mathbf{\nu}$ is the L\'evy measure of an $\alpha$-stable subordinator $\mathbf{\nu}(\dd t)=\alpha t^{-\alpha-1}/\Gamma(1-\alpha)$
and  $\mathbf{\Lambda}(\dd \mathbf{u})=\Prob(\mathbf{V}_1\in\dd \mathbf{u})$, where the sequence $\{\mathbf{V}_i\}_{i\geq 1}$ is defined in point b) of Def. \ref{def:lwolw}.
\end{theorem}
\proof
For a proof of Theorem \ref{th:fconv} we refer to the first part of Appendices B and C in \cite{Magdziarz2015}. The results follows immediately when we put $\gamma=1$ (using notation therein).
\proofend

The above presented result needs an extended comment. First of all, the scaling limit processes for LW and OLW, although both having the form of subordination, differ significantly. In the case of LW, the limit process is a right-continuous version of the left-continuous $\alpha$-stable process $\mathbf{L}^{-}_{\alpha}(s)$ subordinated to the inverse $\alpha$-subordinator $S^{-1}_\alpha(t)$. For the OLW scaling limit we obtain the right-continuous $\alpha$-stable process $\mathbf{L}_{\alpha}(s)$ subordinated to the inverse $\alpha$-subordinator $S^{-1}_\alpha(t)$. In both scenarios the inverse $\alpha$-stable subordinator $S_\alpha^{-1}(t)$ plays the role of internal operational time of the system (see \cite{HenryStraka, Meerschaert31, Magdziarz7, Magdziarz8} for similar results). We emphasize that the difference between scaling limits of LW and OLW processes is due to the fact that the LW process is composed of 'wait-jump' events, while OLW - from 'jump-wait' events. For an extended discussion on this matter see \cite{Magdziarz2015, Magdziarz2012, Teuerle2012, TeuerleOver2012, JurlMeer}.

The $\alpha$-stable process $\mathbf{L}_{\alpha}(s)$ occurs as scaling limit of the sum of jumps $\{\mathbf{J}\}_{i\geq 1}$ (see Def. \ref{def:lwolw}). Since the jumps of the underlying LW and OLW are heavy-tailed with parameter $\alpha$, see (\ref{heavytailed}) and (\ref{jumps}), $\mathbf{L}_{\alpha}(s)$ belongs to the domain of attraction of $\alpha$-stable laws. Moreover, the spatial properties of $\mathbf{L}_{\alpha}(s)$ are also inherited from the sequence $\{\mathbf{V}_i\}_{i\geq 1}$, (\ref{jumps}). It appears that the distribution $\mathbf{\Lambda}$ of unit vector $\mathbf{V}_i$, which governs the possible directions of jump $\mathbf{J}_i$ in (\ref{jumps}), controls the possible $d$-dimensional directions of the jumps of $\mathbf{L}_\alpha(s)$.

Moreover, the LW and OLW scaling limits are composed of two process $\mathbf{L}_\alpha(t)$ and $S_\alpha(t)$, which according the their joint L\'evy measure \eqref{jointlevymeasure1} are strongly dependent. The L\'evy measure can be considered as the intensity of jumps of the joint L\'evy process $\left(\mathbf{L}_\alpha(t),S_\alpha(t)\right)$. Therefore we can conclude that the processes $\mathbf{L}_\alpha(t)$ and $S_\alpha(t)$ have jumps of the same length, which occur in the same epochs of time. As mentioned before, the direction of $d$-dimensional jumps of $\mathbf{L}_\alpha(t)$ is controlled by distribution of $\mathbf{V}_i$.

\subsection{Fractional dynamics of the LW and OLW scaling limits including fractional material derivative}

In what follows, we will use the results of the previous sections to establish a link between the scaling limits of LWs and the dynamics driven by fractional material derivative.

It follows from results presented in \cite{Magdziarz2015,JurlMeer} that the joint L\'evy process $\left(\mathbf{L}_\alpha(t),S_\alpha(t)\right)$ has the following Fourier-Laplace exponent
\begin{equation}
\label{FLT}
\psi(\mathbf{k},s)=
\int\limits_{\mathbf{u}\in\mathbb{S}^{d-1}}\left(s-i\langle\mathbf{k}, \mathbf{u}\rangle\right)^\alpha\mathbf{\Lambda}(\dd \mathbf{u}).
\end{equation}
The Fourier-Laplace transform presented in Eq. (\ref{FLT}) uniquely determines the distribution of the scaling limits obtained in Theorem \ref{th:fconv}. Moreover, it appears that it is closely related to the $d$-dimensional fractional material derivative, which is a $d$-dimensional generalization of one-dimensional fractional derivative introduced in \cite{SokolovMetzler}. The $d$-dimensional fractional material derivative is defined as
\begin{equation}
\label{fmd}
\mathbb{D}^{\alpha,\Lambda}_{\mathbf{x},t} p(\mathbf{x},t) =\int\limits_{\mathbf{u}\in\mathbb{S}^{d-1}}\left(\frac{\partial}{\partial t}+\langle\nabla, \mathbf{u}\rangle\right)^\alpha p(\mathbf{x},t) \mathbf{\Lambda}(\dd \mathbf{u})
\end{equation}
for some density function $p(\mathbf{x},t)$. According to (\ref{FLT}), in the Fourier-Laplace space it is equal to
\begin{equation*}
\mathcal{F}_\mathbf{x}\mathcal{L}_t\{\mathbb{D}^{\alpha,\Lambda}_{\mathbf{x},t}p(\mathbf{x},t)\}=\int\limits_{\mathbf{u}\in\mathbb{S}^{d-1}}\left(s-i\langle \mathbf{k}, \mathbf{u}\rangle\right)^\alpha\mathbf{\Lambda}(\dd \mathbf{u})p(\mathbf{k},s).
\end{equation*}
In the above formula $p(\mathbf{k},s)$ is the Fourier-Laplace transform of $p(\mathbf{x},t)$. Moreover, the operator $\nabla$ denotes the gradient $\nabla=\left(\frac{\partial}{\partial{x_1}},\frac{\partial}{\partial{x_2}},\ldots,\frac{\partial}{\partial{x_d}}\right)^T$, therefore $\langle \nabla,\mathbf{u} \rangle$ can be interpreted as the multidimensional directional derivative, see Section 6.5 in \cite{Meerschaert_book} for details and examples.
The integration in (\ref{fmd}) is over $(d-1)$-dimensional sphere $\mathbb{S}^{d-1}$ as $\mathbf{\Lambda}$ is the distribution on $\mathbb{S}^{d-1}$.

As a consequence the densities of the scaling limits of LW and OLW in Theorem \ref{th:fconv} satisfy in the following fractional diffusion equations \cite{Magdziarz2015}, respectively:
\begin{equation}	
\label{felw}
\mathbb{D}^{\alpha,\Lambda}_{\mathbf{x},t}\,p_{1}(\mathbf{x},t)=\delta_\mathbf{0}(\mathbf{x})\frac{t^{-\alpha}}{\Gamma(1-\alpha)},
\end{equation}
\begin{equation}	
\label{feolw}
\mathbb{D}^{\alpha,\Lambda}_{\mathbf{x},t}\, p_{2}(\mathbf{x},t)=\frac{\alpha}{\Gamma(1-\alpha)}\int\limits_t^\infty\int\limits_{\mathbf{u}\in\mathbb{S}^{d-1}}\delta_{\mathbf{u} s}(\dd \mathbf{x})s^{-\alpha} \Lambda(\dd \mathbf{u}) \dd s.
\end{equation}
Solutions of these equations are in a weak sense, meaning that thy are satisfied in the Fourier-Laplace space \cite{JurlMeer}.

In the next section we will introduce the stochastic scheme which leads to the fractional dynamics given by equations similar to (\ref{felw}) and (\ref{feolw}), but with the operator $\mathbb{D}^{\alpha,\Lambda}_{\mathbf{x},t}$ replaced by the fractional material derivative of distributed-order type.

\section{Stochastic foundations of fractional dynamics with fractional material derivative of distributed-order type}\label{sec:3}
\setcounter{section}{3}
\setcounter{equation}{0}\setcounter{theorem}{0}

\subsection{Definitons of generalized L\'evy walk and overshooting L\'evy walk}
\label{sec:3.1}
Before we present our main results of this section, we shortly recall results concerning the stochastic scheme underlying distributed order derivatives.

Let us consider an IID sequence $\{B_i\}_{i\geq 1}$ of random variables, such that $B_1$ is distributed on interval [0,1] according to probability density function $p(\beta)$, which is regularly varying at zero with exponent $\gamma-1$ (here $\gamma>0$), i.e.
\begin{equation}
\label{B_RV}
\lim_{t\to 0} \frac{p(\lambda t)}{p(t)}=\lambda^{\gamma-1} .
\end{equation}
We also assume that the following property holds
\begin{equation}
\label{B_density}
\int_0^1\frac{p(\beta)}{1-\beta}\mathrm{d}\beta<\infty.
\end{equation}
Then, for any $n\geq 1$ let $\{T^{(n)}_i\}_{i\geq 1}$ be an IID sequence of nonnegative random variables, which conditionally on $B_1=\beta$ has the following distribution
\begin{eqnarray}
\Prob(T^{(n)}_{i}>t|B_i=\beta)=\left\{\begin{array}{lll}
1             &\mathrm{for}& 0\leq t<n^{-1/\beta}, \\
n^{-1}t^{-\beta}    &\mathrm{for}& t\geq n^{-1/\beta}.\\
\end{array}\right.
\label{waiting_times_do}
\end{eqnarray}
%
The asymptotic properties of random variables $\left\{T^{(n)}_i\right\}_{i\geq 1}$ have been investigated in \cite{Meerschaert_ultraslow}, we summarize them in the following lemma.
\begin{lemma}
\label{lemma:B_fconv}
For any $n$ let $\left\{T^{(n)}_i\right\}_{i\geq 1}$ be an IID sequence such that condition (\ref{waiting_times_do}) holds. Then the following convergence holds in $\mathbb{J}_1$ topology as $n\to\infty$:
\begin{equation}
\label{B_fconv}
\sum_{i=1}^{\lfloor n t \rfloor}T^{(n)}_i\fconv S^*_\beta (t),
\end{equation}
where the process is a subordinator such that its L\'evy triplet is given by
\begin{equation}
\label{B_levytriplet}
\left[\ \int\limits_{{x}> 0}\frac{{x}}{1+{x}^2}{\nu}_{S^{*}_\beta}(\dd {x}),{0},{\nu}_{S^{*}_\beta}\right],
\end{equation}
with the L\'evy measure $\mathbf{\nu}_{S^*_\beta}$ given by
\begin{equation}
\label{B_levymeasure}
\mathbf{\nu}_{S^*_\beta}(\dd t)=\int\limits_{0}^1 \beta t^{-\beta-1}\, p(\beta)\dd \beta \dd t \quad\mathrm{for}\ t>0 .
\end{equation}
\end{lemma}
\proof
The proofs follows immediately from Theorem 3.4 and Corollary 3.5 of \cite{Meerschaert_ultraslow}.
\proofend



The specific construction of sequence $\{T^{(n)}_i\}_{i\geq 1}$ leads to the fractional dynamics given by distributed-order derivative \cite{Meerschaert_ultraslow}. Therefore, we will use it to construct the appropriate modifications of LW and OLW processes in order to obtain the dynamics governed by fractional material derivative of distributed-order type.

\begin{definition}
\label{def:lwolw_do}
For any $n\geq 1$ let $\left\{\left(\mathbf{J}^{(n)}_i,T^{(n)}_i\right)\right\}_{i\geq 1}$ be a sequence of independent and identically distributed (IID) random vectors such that all the following conditions hold:
\begin{itemize}
\item[a)] $T^{(n)}_i$'s fulfill condition (\ref{waiting_times_do}),
\item[b)] $\Prob(\mathbf{J}^{(n)}_i\in\mathbb{R}^d)=1$ and for each $\mathbf{J}^{(n)}_i$ we put
\begin{equation}
\mathbf{J}^{(n)}_i=\mathbf{V}_iT^{(n)}_i,
\label{jumps_do}
\end{equation}
where $\{\mathbf{V}_i\}_{i\geq 1}$ is an IID sequence of non-degenerate unit vectors from $\mathbb{R}^d$, sequences $\{\mathbf{V}_i\}_{i\geq 1}$ and $\{T^{(n)}_i\}_{i\geq 1}$ are mutually independent.
\end{itemize}
For the renewal process $N^{(n)}_t=\max\{i\in\mathbb{N}_0:T^{(n)}_1+T^{(n)}_2+\ldots+T^{(n)}_i\leq t\}$, we define the following processes:
\begin{equation}
\mathbf{R}^{(n)}(t)=\sum_{i=1}^{N^{(n)}_t}\mathbf{J}^{(n)}_i,
\label{lw_do}
\end{equation}
\begin{equation}
\mathbf{\tilde{R}}^{(n)}(t)=\sum_{i=1}^{N^{(n)}_t+1}\mathbf{J}^{(n)}_i,
\label{olw_do}
\end{equation}
which are called here generalized L\'evy walk (GLW) and generalized overshooting L\'evy walk (GOLW), respectively.
\end{definition}

\subsection{Asymptotic properties of generalized L\'evy walks and  generalized overshooting L\'evy walks}

In the next part we present the asymptotic behaviour of the GLW and GOLW processes introduced in Sec. \ref{sec:3.1}. We give detailed description of the resulting scaling limit processes.

Let us start with the technical lemma concerning the asymptotic properties of the joint process of partial sums for $\{(\mathbf{J}^{(n)}_i,T^{(n)}_i)\}_{i\geq 1}$.
\begin{lemma}
\label{lemma:joint_fconv}
Let us consider an IID sequence $\{(\mathbf{J}^{(n)}_i,T^{(n)}_i)\}_{i\geq 1}$ from Def. \ref{def:lwolw_do}. Then the following convergence holds in $\mathbb{J}_1$ topology as $n\to\infty$:\begin{equation}
\label{joint_fconv}
\left(\sum_{i=1}^{\lfloor n t \rfloor}\mathbf{J}^{(n)}_i,\sum_{i=1}^{\lfloor n t \rfloor}T^{(n)}_i\right)\fconv \left(\mathbf{L}^*_\beta (t),S^*_\beta (t)\right),
\end{equation}
where the process $\left(\mathbf{L}^*_\beta(t),S^*_\beta (t)\right)$ is (d+1)-dimensional L\'evy process defined by the following L\'evy triplet:
\begin{equation}
\label{joint_levytriplet}
\left[\ \int\limits_{\mathbf{x}\neq \mathbf{0}}\frac{\mathbf{x}}{1+\|\mathbf{x}\|^2}{\nu}_{\left(\mathbf{L}^{*}_\beta,S^{*}_\beta\right)}(\dd {x}),{0},{\nu}_{\left(\mathbf{L}^{*}_\beta,S^{*}_\beta\right)}\right],
\end{equation}
with the following L\'evy measure
\begin{equation}
\label{joint_levymeasure}
\mathbf{\nu}_{\left(\mathbf{L}^{*}_\beta,S^{*}_\beta\right)}(\dd \mathbf{x})=\int\limits_{0}^1 \left(\int_{\mathbf{u}\in \mathbb{S}^{d-1}}\delta_{t\mathbf{u}}(\dd \mathbf x_1) \mathbf{\Lambda}(\dd \mathbf{u}) \beta t^{-\beta-1}\dd t\,  \right)p(\beta)\dd \beta.
\end{equation}
Here $\mathbf{x}=(\mathbf{x_1},t)\in \mathbb{R}^d\setminus\{\mathbf{0}\}\times(0,\infty)$, $\mathbf{u}=\mathbf{x_1}/\|\mathbf{x_1}\|\in \mathbb{S}^{d-1}$, $\mathbb{S}^{d-1}$ is (d-1)-dimensional sphere, $\delta_{s}(\cdot)$ is a Dirac delta function at point $s$ and  $\mathbf{\Lambda}(\dd \mathbf{u})=\Prob(\mathbf{V}_1\in\dd \mathbf{u})$, where the sequence $\{\mathbf{V}_i\}_{i\geq 1}$ is defined in point b) of Def. \ref{def:lwolw_do}.
\end{lemma}
\proof
First, we recall a result concerning finite-dimensional convergence of waiting times defined in (\ref{waiting_times_do}). Namely, from Lemma \ref{lemma:B_fconv} it follows that
$$\sum_{i=1}^{ n } T^{(n)}_i \stackrel{d}{\longrightarrow} S^*_\beta(1),$$
holds as $n\to\infty$. Here $S^*_\beta(1)$ is the distribution of the process $S^*_\beta(t)$ at $t=1$. Process $S^*_\beta(t)$ is fully described by its L\'evy triplet given by (\ref{B_levytriplet}) with the L\'evy measure ${\nu}_{S^{*}_\beta}$ given by (\ref{B_levymeasure}).
Based on Theorem 3.2.2 in \cite{Meerschaert1} the above mention result implies that as $n\to\infty$
\begin{equation}
\label{conv_waitingtimes_do}
n \Prob(T_1^{(n)}\in\dd t)\to {\nu}_{S^{*}_\beta}(\dd t)=\int_0^1\beta t^{-\beta-1} p(\beta) \dd\beta \dd t \quad\mathrm{for}\ \mathrm{any}\ t>0
\end{equation}
and the Gaussian part in L\'evy-Khinchin representation of the process $S^*_\beta(t)$ is equal to 0.

To prove Lemma \ref{lemma:joint_fconv}, let us observe that for any Borel sets $\mathbf{A}_1\in\mathcal{B}(\mathbb{R}^d)$ and $A_2 \in \mathcal{B}(\mathbb{R}_+ )$ such that  $\mathbf{A}_1=\mathbf{A}_1(R,\mathbf{D})=\{r\mathbf{u}\in\mathbb{R}^d:r\in R, \mathbf{u} \in \mathbf{D}\}$, where $ R\in\mathcal{B}(\mathbb{R}_+)$ and $\mathbf{D}\in\mathcal{B}(\mathbb{S}^{d-1})$, due to the independence between $T^{(n)}_1$ and $\mathbf{V}_1$ we have that
\begin{eqnarray}
n\Prob(\mathbf{J}^{(n)}_1\in \mathbf{A}_1, T^{(n)}_1 \in A_2)\nonumber
=n\! \int_{\mathbb{S}^{d-1}}\!\! \Prob(T^{(n)}_1 \mathbf{u}\in \mathbf{A}_1, T^{(n)}_1 \in A_2)\Prob(\mathbf{V}_1\in \dd \mathbf{u}) \\\nonumber
=n\! \!\!\int_{A_2}\int_{\mathbb{S}^{d-1}}\!\! \indyk(t\mathbf{u}\in\mathbf{A}_1)\Prob(\mathbf{V}_1\in \dd \mathbf{u})\Prob(T^{(n)}_1\in \dd t) \\ \label{joint_conv_calc1}
\end{eqnarray}
where $\indyk({\cdot})$ is equal 1 if the condition in the brackets is fulfilled and 0 if it is not.
Let us observe that for $(\mathbf{x},t)\in\mathbb{R}^d\setminus\{\mathbf{0}\}\times\mathbf{R}_+$ and $\mathbf{u}\in\mathbb{S}^{d-1}$ by taking the advantage of Fubini's theorem and asymptotic property (\ref{conv_waitingtimes_do}) it follows from the formula (\ref{joint_conv_calc1}) that
\begin{eqnarray}
n\Prob(\mathbf{J}^{(n)}_1\in \mathbf{A}_1, T^{(n)}_1 \in A_2)=n\! \int_{A_2}\int_{\mathbb{S}^{d-1}}\!\! \indyk(t\mathbf{u}\in\mathbf{A}_1)\Prob(\mathbf{V}_1\in \dd \mathbf{u})\Prob(T^{(n)}_1\in \dd t)\nonumber\\\nonumber
\longrightarrow \int_0^1 \int_{A_2} \int_{\mathbb{S}^{d-1}}\indyk(t\mathbf{u}\in\mathbf{A}_1) \beta t^{-\beta-1}\mathbf{\Lambda}(\dd \mathbf{u})p(\beta)\dd t \dd\beta,\\\label{joint_conv_calc2}
\end{eqnarray}
as $n\to\infty$, where $\mathbf{\Lambda}(\dd \mathbf{u})=\Prob(\mathbf{V}_1\in \dd \mathbf{u})$.
Thus the L\'evy measure of joint distribution $(\mathbf{L}^*_\beta(1),S^*_\beta (1))$ is equal to 
\begin{equation}
\nonumber
\mathbf{\nu}_{\left(\mathbf{L}^{*}_\beta,S^{*}_\beta\right)}(\dd \mathbf{x})=\int\limits_{0}^1 \left(\int_{\mathbf{u}\in \mathbb{S}^{d-1}}\delta_{t\mathbf{u}}(\dd \mathbf x_1) \mathbf{\Lambda}(\dd \mathbf{u}) \beta t^{-\beta-1}\dd t\,  \right)p(\beta)\dd \beta.
\end{equation}
Here $\mathbf{x}=(\mathbf{x_1},t)\in \mathbb{R}^d\setminus\{\mathbf{0}\}\times(0,\infty)$, $\mathbf{u}=\mathbf{x_1}/\|\mathbf{x_1}\|\in \mathbb{S}^{d-1}$.

Since the distribution $S^*_\beta (1)$ has Gaussian component equal to 0, one can conclude that due to the tight relation between ${\mathbf{J}^{(n)}_i}$ and ${{T}^{(n)}_i}$, also the joint limit distribution $(\mathbf{L}^*_\beta(1),S^*_\beta (1))$ has no Gaussian component.
Theorem 3.2.2 \cite{Meerschaert1} assures that the convergence in (\ref{joint_fconv}) holds in the sense of distributions. Based on Theorem 4.1 in \cite{MeerschaertInvariance} we conclude that the convergence is also valid in the functional sense (here convergence in Skorokhod $\mathbb{J}_1$ topology, see \cite{Billingsley,Whitt}), which finishes the proof.\proofend

In the next theorem we present the main result of this section, which are the scaling limits of GLW and GOLW processes.
\begin{theorem}
\label{th:fconv_do}
Let $\mathbf{R}^{(n)}(t)$ and $\mathbf{\tilde{R}}^{(n)}(t)$ be the GLW and GOLW processes defined in Def. \ref{def:lwolw_do}. The following convergences hold in $\mathbb{J}_1$ topology as $n\to\infty$:
\begin{eqnarray}
\label{th:fconv1_do}
{\mathbf{R}^{(n)}(t)}\Rightarrow \left( \mathbf{L}_\beta^{*-}\left( S^{*-1}_{\beta}(t)\right)\right)^+,\\\nonumber {\mathbf{\tilde{R}^{(n)}}(t)}\Rightarrow  \mathbf{L}^{*}_{\beta}\left( S^{*-1}_{\beta}(t)\right),
\end{eqnarray}
where the process ${S}^{*-1}_{\beta}(t)$ is defined as follows:
\begin{equation}
{S}^{*-1}_{\beta}(t)=\inf\{\tau \geq 0: {S}^*_{\beta}(\tau)>t\}.
\label{inversestablesubordinator_do}
\end{equation}
The dependence structure between $\mathbf{L}^{*}_{\beta}(t)$ and $S^{*}_{\beta}(t)$ processes is given by their joint L\'evy triplet in (\ref{joint_levytriplet}).
\end{theorem}
\proof
In Lemma \ref{lemma:joint_fconv} we showed the functional convergence of the joint process of cumulated jumps and waiting times to $(L^*_\beta(t),S^*_\beta(t))$, which is fully described by its L\'evy triplet given by (\ref{joint_levytriplet}). Moreover, from the L\'evy triplet of $S^*_\beta(t)$ we can conclude that it is a process with unbounded, nonnegative and strictly increasing trajectories. Therefore, the convergence in Skorokhod $\mathbf{J}_1$ topology in Theorem \ref{th:fconv_do} can be now easily concluded from Theorem 3.6  in \cite{HenryStraka}.
\proofend

The general structure of the scaling limit processes in Theorem \ref{th:fconv_do} is similar to the one in Theorem \ref{th:fconv}. We see the difference for GLW and GOLW scaling limits due to the fact that these processes arise form different stochastic schemes: 'wait-first' and 'jump-first' scenarios, respectively. For the GLW scheme L\'evy process $L^{*-}_{\beta}(s)$ is subordinated to the inverse subordinator $S^{*-1}_\beta(t)$. While for the GOLW scaling limit process we obtain right-continuous L\'evy process $L^*_{\beta}(s)$ subordinated to the inverse subordinator $S^{-1}_\beta(t)$. Similarly as for $S_\alpha^{-1}(t)$ in Theorem \ref{th:fconv}, here the inverse subordinator $S_\beta^{*-1}(t)$ plays a role of internal operational time of the resulting scaling limit processes.

It is worth noticing that processes $L^*_{\beta}(s)$ and $S_\beta^{*}(t)$ are no longer stable processes. Due to the specific construction of waiting times, see (\ref{waiting_times_do}) and (\ref{heavytailed}), where the heavy-tailed index $\beta$  being conditionally defined by the sequence $\{B_i\}_{i\geq 1}$, $L^*_{\beta}(s)$ and $S_\beta^{*}(t)$ belong to the more general class of L\'evy processes. Processes $L^*_{\beta}(s)$ and $S_\beta^{*}(t)$ are strongly dependent, which can be deduced from their joint L\'evy triplet. It appears  that their joint L\'evy measure (\ref{joint_levymeasure}) is an integral over the distribution of parameter $\beta$ of joint L\'evy measures (\ref{jointlevymeasure1}) in the LW and OLW case

\subsection{Fractional dynamics of the GLW and GOLW scaling limits including fractional material derivative of distributed-order type}
It follows from results presented in \cite{Magdziarz2015,JurlMeer} that the scaling limit of processes $\left(\mathbf{L}^*_\beta(t),S^*_\beta(t)\right)$ has the following Fourier-Laplace exponent
\begin{equation}
\label{FLT_do}
\psi(\mathbf{k},s)=
\int_0^1\left(\int\limits_{\mathbf{u}\in\mathbb{S}^{d-1}}\left(s-i\langle\mathbf{k}, \mathbf{u}\rangle\right)^\beta \Gamma(1-\beta)\mathbf{\Lambda}(\dd \mathbf{u})\right)p(\beta)\dd \beta.
\end{equation}
The Fourier-Laplace transform presented in Eq. (\ref{FLT_do}) uniquely determines distributions of the corresponding scaling limits obtained in Theorem \ref{th:fconv_do}. Therefore, one can introduce a pseudodifferential operator of fractional material derivative of distributed-order type:
\begin{equation}
\label{fmd_do}
\mathbb{D}^{p(\beta),\Lambda}_{\mathbf{x},t} p(\mathbf{x},t)=\int_0^1\left(\int\limits_{\mathbf{u}\in\mathbb{S}^{d-1}}\left(\frac{\partial}{\partial t}+\langle\nabla, \mathbf{u}\rangle\right)^\beta\Gamma(1-\beta) p(\mathbf{x},t)\mathbf{\Lambda}(\dd \mathbf{u})\right)p(\beta)\dd \beta
\end{equation}
which according to (\ref{FLT_do}) is defined for some density function $p(\mathbf{x},t)$ in the Fourier-Laplace space as
\begin{equation*}
\mathcal{F}_\mathbf{x}\mathcal{L}_t\{\mathbb{D}^{p(\beta),\Lambda}_{\mathbf{x},t}p(\mathbf{x},t)\}=\int_0^1\left(\int\limits_{\mathbf{u}\in\mathbb{S}^{d-1}}\left(s-i\langle \mathbf{k}, \mathbf{u}\rangle\right)^\beta\Gamma(1-\beta)\mathbf{\Lambda}(\dd \mathbf{u})\right)p(\mathbf{k},s) p(\beta)\dd \beta .
\end{equation*}
In formula (\ref{fmd_do}) the first integral is the one with respect to the distribution $p(\beta)$ over the integrand which is the usual fractional material derivative operator (\ref{fmd}).

In analogy to the result of Section \ref{sec:2} the densities of the scaling limits of GLW and GOLW in Theorem \ref{th:fconv_do} satisfy the following fractional diffusion equations, respectively:
\begin{equation}	
\label{felw_do}
\mathbb{D}^{p(\beta),\Lambda}_{\mathbf{x},t}\,p_{1}(\mathbf{x},t)=\delta_\mathbf{0}(\mathbf{x})\int_0^1{t^{-\beta}}p(\beta)\dd \beta,
\end{equation}
\begin{equation}	
\label{feolw_do}
\mathbb{D}^{p(\beta),\Lambda}_{\mathbf{x},t}\, p_{2}(\mathbf{x},t)=\int_0^1\int_t^\infty\int_{\mathbf{u}\in\mathbb{S}^{d-1}}\delta_{s \mathbf{u}}(\mathbf{x}) s^{-\beta}\Lambda(\dd \mathbf{u}) \dd s p(\beta) \dd \beta,
\end{equation}
where the fractional material derivative operator of distributed order is on the left-hand side of above equations.

\section{Conclusions}
In this paper we introduced the Levy walk model of distributed-order type.
Our approach was based on strongly coupled CTRW with the distribution of waiting times
displaying ultraslow (logarithmic) decay of the tails. This type of distribution
is closely related to the so-called ultraslow diffusion.
We have derived the diffusion limit of the Levy walk model of distributed-order type.
It appears that the limit has the form of subordination (time change) of certain L\'evy process
with the inverse subordinator corresponding to the waiting times with ultra-heavy tails.
This type of superdiffusive dynamics is characterized by trajectories that have very long jumps (with logarithmic decay of their tails)
and finite mean square displacement (in the case of GLW).

The introduced model explains the stochastic origins of fractional dynamics driven by fractional material derivative of distributed order-type.
This is manifested by the fact that the probability density function of the obtained diffusion limit process solves the fractional diffusion equation with such material derivative.

 \medskip

\section*{Acknowledgements}

This research was partially supported by NCN Maestro grant no. 2012/06/A/ST1/00258.



 \bigskip \smallskip

 \it

 \noindent
Hugo, Steinhaus Center, Department of Mathematics \\
Wroclaw University of Technology \\
Wybrze\.ze Wyspia\'nskiego 27 \\
50-370 Wroc\l aw, POLAND  \\[4pt]
$^1$  e-mail: marcin.magdziarz@pwr.edu.pl
$^2$  e-mail: marek.teuerle@pwr.edu.pl

\end{document}